\documentclass{article}

\usepackage{amssymb}
\usepackage{amsmath}

\usepackage[super]{nth}

\usepackage{sectsty}
\sectionfont{\fontsize{15.8}{14}\selectfont}

\usepackage{cite}

\usepackage{hyperref}
\hypersetup{
    colorlinks=true,
    citecolor=black,
    linkcolor=black
}

\usepackage{graphicx}
\graphicspath{ {./images/} }

\newcommand{\La}{\lambda}
\newcommand{\sLa}{{\scriptscriptstyle\lambda}}
\newcommand{\Na}{\mathbb{N}}
\newcommand{\QQ}{\mathbb{Q}}
\newcommand{\BB}{\mathbb{B}}
\newcommand{\RR}{\mathbb{R}}
\newcommand{\ZZ}{\mathbb{Z}}
\newcommand{\bb}{\mathcal{B}}
\newcommand{\ud}{\updownarrow}

\newcommand{\bo}{\boldsymbol}

\usepackage{tgadventor}

\begin{document}
\title{At the End of Infinity}
\author{Emmanuel Rochette}

\begin{center}
\LARGE\textbf{At the End of Infinity} \\ \vspace{2.23mm}
\large{Emmanuel Rochette} \\
\small{March 8, 2021} \\ \vspace{1.7mm}
\small{emmanuel.rochette@mail.mcgill.ca} \vspace{1.3mm}
\end{center}

\noindent
Mathematics is full of paradoxes and limitations.  As this work will demonstrate, many of such problems arise from a general misunderstanding of how to properly use infinity. In particular, we'll argue that set theory, which relies heavily on Georg Cantor's ideas, is simply inconsistent with Calculus \textbf{---} and, to remedy this situation, we'll develop our own theory of infinity. 

To discuss the properties of an infinite set, mathematicians currently rely on two kinds of numbers \cite{settheory}. They use the \textit{cardinals} to count how many objects there are in a set, with the landmark result that {\small$\Na$}'s cardinality is smaller than {\small$\RR$}'s. In all cases, however, an infinite quantity is believed to remain unchanged after any arithmetical operation, and they thus commonly use \vspace{0.5mm}
\begin{equation*}
a\cdot\infty^b + c = \infty
\end{equation*}
The \textit{ordinals} are also used to say what position each object has in a set, relative to one another \textbf{---} \nth{1}, \nth{2}, \nth{3}, etc. Mathematicians find the ordinals useful because they can enumerate objects past infinity; after the infinite ordinal {\small$\omega$}, a new enumeration can begin, and so on, as to give the following sequence \vspace{0.4mm}
\begin{equation*}
1, \: 2  \: ,3,  \: \ldots,  \: \omega,  \: \omega+1, \: \omega+2, \: \ldots, \: \omega2, \: \ldots,  \: \omega n,  \: \ldots,  \: \omega^2 ,  \: \ldots,  \: \omega^\omega,  \: \ldots
\end{equation*}
But it's assumed that placing new objects in front of others will not change the enumeration process, so the ordinals are not commutative in both addition and multiplication. More specifically, they write \vspace{0.4mm}
\begin{equation*}
n + \omega = \omega \quad\quad \text{and} \quad\quad \omega + n \neq \omega 
\end{equation*}

The distinction between cardinals and ordinals came from Cantor's work, and it therefore has almost 150 years of mathematical history \cite{mathhistory, CantorOrigin}. We'll demonstrate, nevertheless, that many severe inconsistencies emerge directly from it \textbf{---} and we'll thus argue against making such a distinction.   

In fact, to resolve these issues, we'll introduce our own infinite number {\small$\La$}, which will behave as both a cardinal and an ordinal number. Its arithmetical properties will be developed in the five upcoming sections, where each of them respectively presents one basic operation of arithmetic.  

Our work will strongly disagree with many well-established results from Cantor's theory; we'll even refute the famous one-to-one correspondence between the sets {\small$\Na$}, {\small$\ZZ$} and {\small$\QQ$}. Since this might provoke controversy, we pledge to follow two important guidelines. First, we'll keep our presentation as simple and accessible as possible. But, more importantly, we'll use Calculus to verify each new property of our number {\small$\La$} \textbf{---} and, when it's applicable, we'll even use it to resolve long-standing paradoxes in mathematics.

\section*{1. Addition  }
The infinite number {\small$\La$} will be the main focus of this work, and we thus start by introducing it. First, consider the set {\small$\Na_n=\{1,2,3,...,n\}$}. Clearly, its cardinality is 
\begin{equation}
\big| \Na_n \big| = n \\[3.5pt]
\end{equation}
Since $n$ is arbitrary, it follows that all subsets of {\small$\Na$} have their cardinality equal to their last number. From this observation, after taking the limit of $n$ goes to infinity, we are lead to write \vspace{0.4mm}
\begin{equation}
\Na_\La = \Na = \{1,2,3,...,\La \} 
\end{equation}
with \vspace{0.3mm}
\begin{equation} 
\big| \Na_\La  \big| = \La \\[4pt]
\end{equation}
The number {\small$\La$} will thus be equal, throughout this work, to the cardinality of {\small$\Na$}. Similarly to the ordinal numbers, to extend it further, we can use the operation of addition and define {\small$\La+1$}, then {\small$\La+2$}, and so on, up to yet another limit to infinity, {\small$\La+\La=2\La$}. The set {\small$\Na_\La$} is thus extended, after taking this second infinite limit, to the larger set  \vspace{0.2mm}
\begin{equation}
\Na_{2\La} = \{1,2,3,...,\La\} \cup \{\La+1, \La+2,...,2\La\}
\end{equation}
with, evidently \vspace{0.4mm}
\begin{equation} 
\big|\Na_{2\La}\big| = \La+\La = 2\La  \\[5.8pt]
\end{equation}

Recognizably, this set is the equivalent of the already known {\small$\omega\cdot2$}. In mathematics, however, the sets {\small$\Na$} and {\small$\omega\cdot2$} are said to have the same cardinality, and no mathematician ever uses the ordinal numbers to measure different sizes of infinities. As we'll demonstrate in this work, this is unfortunately due to many wrong results from Cantor's theory, where the one-to-one correspondences are mainly to blame. As a first example, consider the two infinite sets \vspace{0.6mm}
\begin{align}
\begin{split} \label{eq:omegamapping}
\Na = \{2,4,6,  . . . \} \cup \{1,3,5, . . . \} \quad\quad\quad\; \\[3.5pt]
\text{and} \quad\quad\quad\quad\quad\quad\quad\quad\;\: \\[3pt]
\omega \cdot 2 = \{1, 2, 3, . . . \} \cup \{\omega, \omega + 1, \omega + 2, . . .\}  \\[1.5pt]
\end{split}
\end{align}
where all {\small$\Na$}'s even numbers would be mapped to {\small$\{1, 2,3, . . . \}$}, and where all odd numbers would have a correspondence with {\small$\{\omega, \omega + 1, \omega + 2, . . .\}$}. 

Because of this mapping in particular, most mathematicians currently believe, paradoxically, that there are as many even numbers in {\small$\Na$} as there are numbers in this same set \textbf{---} although some are even while others are odd. That's an evident self-contradiction, but it's usually justified by simply saying that infinity has a counterintuitive nature. Nevertheless, we'll now demonstrate, using Calculus, that such a paradox just cannot be true. 

To do so, we'll rely on Riemann's definition of an integral. Typically, it's given as \vspace{0.4mm}
\begin{equation}
\int_{a}^{b} f(x)dx = \lim_{n\to\infty} \sum_{j=1}^{n} f(a+jdx) \cdot dx \\[2.5pt]
\end{equation}
where $dx=(b-a)/n$. In fact, we want to compare two integrals, the first one with $f_1(x) = x$, and the other with $f_2(x) = 2x$. From the above equation, these are written as
\begin{align}
\begin{split} \label{eq:1}
\quad \int_{0}^{b} f_1(x)dx &=  \lim_{n\to\infty} \sum_{j=1}^{n} (jdx) \cdot dx \\[3pt]
			       &= (1 + 2+ 3 + 4 +5+6+\cdots ) \cdot dx^2 \\[6.5pt]
			       &=(1+3+5+\cdots) \cdot dx^2  \\
			       & \quad\quad\quad + (2+4+6+\cdots) \cdot dx^2  \\[4.8pt]
			       &= b^2/2
\end{split}
\end{align}
and 
\begin{align}
\begin{split} \label{eq:2}
\int_{0}^{b} f_2(x)dx &= \lim_{n\to\infty} \sum_{j=1}^{n}(2jdx) \cdot dx \\[2.5pt]
			       &= (2+4+6+\cdots) \cdot dx^2 \\[6pt]
			       &= b^2
\end{split}
\end{align}
where the infinite limits were taken, just for the moment, without any reference to {\small$\La$}. Also, the Riemann sum in 1.\ref{eq:1} was rearranged in a suggestive way. 

Now, for any {\small$b>0$}, notice that the integral with {\small$f_1(x)$} gives a smaller value than the one with {\small$f_2(x)$}. Their dependance on {\small$b$} comes, however, from the same factor {\small$dx^2=(b/n)^2$}. Consequently, their difference in value has to come from their respective infinite sum of numbers.

We can now make two important observations. First, since the index $j$ runs from $1$ to infinity, the Riemann sum in 1.\ref{eq:1} must contain all of {\small$\Na$}'s numbers. Secondly, according to Cantor's theory, there's a one-to-one correspondence between {\small$\Na$} and its even numbers; the Riemann sum in 1.\ref{eq:2} would therefore only contain those.

A very unfortunate problem now arises, however: the Riemann sum in 1.\ref{eq:1} contains all the even and odd numbers of {\small$\Na$}, but it would still give a smaller value than 1.\ref{eq:2}'s, which only contains the even numbers. 

Nevertheless, if we re-introduce {\small$\La$}, this would-be paradox has a very simple resolution. In particular, the set {\small$\Na$} does not contain all the terms in 1.\ref{eq:2} \textbf{---}  it's the larger set {\small$\Na_{2\La}$} that actually contains all of them. To see this, we must stop using the {\small$+\cdots$} notation, and work with {\small$\La$} instead. This allows us to write \vspace{0.01mm}
\begin{align}
\begin{split} \label{eq:3}
\int_{0}^{b} f_1(x)dx &=  \sum_{j=1}^{\La} (jdx) \cdot dx  \\[4pt]
			       &= (1 + 2+ 3 + 4 + \cdots + \La) \cdot dx^2 \\[5.5pt]
			       &= b^2/2
\end{split}
\end{align}
and 
\begin{align}
\begin{split} \label{eq:4}
\int_{0}^{b} f_2(x)dx &= \sum_{j=1}^{\La} (2jdx) \cdot dx \\[4pt]
			       &= (2+4+6+ 8 +\cdots+2\La) \cdot dx^2 \\[5.5pt]
			       &= b^2 
\end{split}
\end{align}
It's already clear that {\small$2\La \in \Na_{2\La}$}, although {\small$2\La \notin \Na_\La$}. In fact, the Riemann sum in 1.\ref{eq:4} contains infinitely many terms that are in the set {\small$\Na_{2\La}$}, but not in {\small$\Na_{\La}$}. While the exact proportions will be given in Section $4$, we can nonetheless explain why the Riemann sum in 1.\ref{eq:4} gives a bigger value than 1.\ref{eq:3}'s: all the numbers {\small$\La+n$} from {\small$\Na_{2\La}$} are bigger than any other from {\small$\Na_{\La}$}, and they thus contribute more to the integral's value.  

Now, for yet another inconsistency between Calculus and Cantor's theory, consider the sets {\small$\Na$} and {\small$\ZZ$}. They are claimed to have the same cardinality, due to an alleged one-to-one correspondence between them. Typically, if {\small$\ZZ$} is taken without $0$ for simplicity, the mapping is done as shown below.  \vspace{0.75mm}
\begin{equation*}
\begin{matrix}
\quad\quad\quad \mathbb{N}    &\bo{:}   &1\;\;     &2\;\;       &3      &\:\:4               &\:\:5         &\:\:6           &  \\[6pt]
\quad\quad\quad  		     &           &\bo{\ud}\;\;   &\bo{\ud}\;\;   &\bo{\ud}   &\:\:\bo{\ud}            &\:\:\bo{\ud}      &\:\:\bo{\ud}         &\bo{\cdots}  \\[7pt]
\quad\quad\quad \ZZ 		     &\bo{:}   &1\;\;     &-1\:\:\;\:      &2    &\:\:-2\:\:       &\:\:3  	  &\:\:-3\:\:\:	    &
\end{matrix}
\end{equation*}
\vspace{0.2mm}But, here again, this is not one-to-one: as we enumerate {\small$\Na$} from 1 to {\small$n$}, let's say, the mapping only reaches the integers from {\small$-n/2$} to {\small$n/2$}. Consequently, there's a total of {\small$n$} integers that are yet to appear in the enumeration of {\small$\ZZ$}. In fact, as we continue to enumerate all of {\small$\Na$}, this gap will never decrease; only {\small$\La$} integers are reached by this mapping, and {\small$\La$} others are left out.

To fix this problem, we simply need to extend {\small$\Na$} to the set {\small$\Na_{2\La}$}. The mapping can thus be \vspace{0.9mm}
\setcounter{MaxMatrixCols}{14}
\begin{equation*}
\begin{matrix}
\quad\: \mathbb{N}_{2\La} &\bo{:}   &1\;\;     &2\;\;      &3       &                     &\La   &\;            &{\La{\scriptstyle+}1}    &\:{\La{\scriptstyle+}2}    &\:{\La {\scriptstyle+}3}      &                     &2\La \\[7pt]
\quad\:  		             &            &\bo{\ud}\;\;  &\bo{\ud}\;\;   &\bo{\ud}    &\bo{\cdots}    &\bo{\ud}   &\bo{,}\;    &\bo{\ud}                              &\:\bo{\ud}                              &\:\bo{\ud}                                  &\bo{\cdots}     & \bo{\ud}  \\[8pt]
\quad\: \ZZ 		             &\bo{:}   &1\;\;     &2\;\;      &3       &                     &\La   &\;            &-1\:\:\:                           &\:-2\:\:\:                         &\:-3\:\:\:                             &                      &-\La\:\:\: 
\end{matrix}
\end{equation*}
\pagebreak

\noindent
From this one-to-one correspondence, the final conclusion must be  \vspace{0.3mm}
\begin{equation}
|\ZZ| = |\Na_{2\La}| = 2|\Na| = 2\La 
\end{equation}

The cardinality of {\small$\ZZ$} is thus twice the size of {\small$\Na$}. In fact, we can demonstrate this result directly from the well-known identity \vspace{0.2mm}
\begin{equation} \label{eq:a6}
\int_{-\infty}^{\infty} f(x) dx = 2\int_{0}^{\infty} f(x) dx 
\end{equation}
where {\small$f(x)$} is any even function. First, we need to partition these integrals into unit intervals. The above equation will thus become \vspace{1mm}
\begin{equation} \label{eq:bigass}
\sum_{n=1}^\La \int_{n-1}^n f(x)dx + \sum_{n=1}^\La \int_{-n}^{-n+1} f(x)dx = 2 \cdot \sum_{n=1}^\La \int_{n-1}^n f(x)dx \\[1.5pt]
\end{equation}
Clearly, there's a direct one-to-one mapping between {\small$\ZZ$} and the integrals of the {\small LHS}: the {\small$n$}th positive integer is mapped to {\small$\int^n_{n-1} f(x)dx$}, and the {\small$n$}th negative integer is mapped to {\small$\int^{-n+1}_{-n} f(x)dx$}. Similarly, there's a one-to-one correspondence between {\small$\Na$} and the integrals in the {\small RHS}. 

However, since {\small$f(x)$} is an even function, it follows that \vspace{0.3mm}
\begin{equation}
\int^{-n+1}_{-n} f(x)dx = \int^{n+1}_{n} f(x)dx
\end{equation}
The equality in 1.\ref{eq:bigass} therefore holds because there are twice as many integrals in the {\small RHS} as there are in the {\small LHS}. Consequently, due to their respective one-to-one correspondence with {\small$\ZZ$} and {\small$\Na$}, we must conclude in the same way for their cardinality: there are twice as many integers as there are natural numbers. 

Now, before concluding this section, there's a last paradox that we want to resolve, since it's very similar to the problematic mapping in 1.\ref{eq:omegamapping}. It was first introduced in 1638 by Galileo, when he published \textit{Two New Sciences}, his last scientific work.  

The reasoning goes as follows. First, since some numbers are squares while others are not, Galileo concludes that {\small$\Na$} contains less squared numbers than there are in total. However, since every number can be squared, it also appears that {\small$\Na$} contains as many numbers as there are squared ones. Galileo's final conclusion was that "the attributes equal, greater, and less, are not applicable to infinite, but only to finite, quantities".

Nevertheless, more than two centuries later, Cantor favoured the alleged one-to-one correspondence between {\small$\Na$} and its squared numbers, which is now accepted by most mathematicians. By using Calculus again, however, we can show that Cantor's claim is incorrect \textbf{---} and that Galileo's paradoxical conclusion is wrong too. 
 
Simply consider the following integral. 
\begin{align}
\begin{split} \label{eq:5}
\int_{0}^{b} x^2dx &= \sum_{j=1}^{\La} (jdx)^2 \cdot dx \\[4pt]
			       &= (1+4+9+16+\cdots+\La^2) \cdot dx^3 \\[5.5pt]
			       &= b^3/3 \\[0.2pt]
\end{split}
\end{align}
Since the Riemann sum in 1.\ref{eq:3} contains all the numbers of {\small$\Na$}, we'll compare it to the above one. First, for {\small$b>3/2$}, notice that {\small$b^2/2<b^3/3$}. But {\small$dx^3<dx^2$} for any {\small$b<\La$}. Consequently, we must already conclude that \vspace{0.3mm}
\begin{equation} 
(1+2+3+4+\cdots+\La) < (1+4+9+16+\cdots+\La^2) \\[0.8pt]
\end{equation}
For this inequality to hold, the RHS must contain terms that belong to a bigger set than {\small$\Na$}, since the LHS already contains all of {\small$\Na$}'s numbers. In particular, we must extend the set {\small$\Na_\La$} to the larger {\small$\Na_{\La^2}$}, which includes the square of all numbers from $1$ to {\small$\La$}. We therefore have a set that resolves Galileo's paradox and refutes Cantor's claim, and the quantity {\small$\La^2 = \La \cdot \La$} now brings us to the next section. \vspace{11mm}  

\setcounter{equation}{0} 
\setcounter{section}{2}

\section*{2. Multiplication}
In the previous section, we started by introducing {\small$\La$} and some of its basic attributes, and we used it to resolve a few inconsistencies in mathematics. In this section, however, we'll proceed in reverse: the currently accepted properties of infinity will be shown, once again, to be incompatible with Calculus, and we'll expand our theory to fix those issues. In particular, we'll need to extend the set {\small$\Na_{2\La}$} even further \textbf{---} and this will be our opportunity to formally discuss how to multiply with {\small$\La$}. 

Now, to start exposing these inconsistencies, we simply have to compare the two integrals below. \vspace{0.6mm}
\begin{equation}
\int_0^b xdx \quad \text{and} \quad \int_0^b \int_0^bdxdy \label{eq:0}  \\[4.5pt]
\end{equation}
In 1.\ref{eq:3}, we gave the Riemann sum of the first integral, and it was very simple: {\small$(1+2+3+\cdots +\La) \cdot dx^2 $}. For the other one, the Riemann sum is a bit more involved, since it's actually made from two infinite sums. To simplify the notation, we'll thus add the subscript {\small$\La$} to the last term of a sum to indicate that it contains {\small$\La$} terms. Our double Riemann sum can thus be given as 
\begin{align}  \label{eq:b1}
\begin{split} 
\int_0^b \int_0^b dxdy &= \int_0^b \Big( (1+1+1+\cdots+1_{\scriptscriptstyle\lambda})  dx \Big)  dy \\[5pt]
				  &=
 		 \begin{pmatrix}
		\:\:\:\:\:(1+1+1+1+\cdots+1_{\scriptscriptstyle\lambda}) \\[1.5pt]
				  +\: (1+1+1+1+\cdots+1_{\scriptscriptstyle\lambda}) \\[1.5pt]
				+\:(1+1+1+1+\cdots+1_{\scriptscriptstyle\lambda})  \\[1.5pt]
				+\:(1+1+1+1+\cdots+1_{\scriptscriptstyle\lambda})  \\[1.5pt]
				   \;\;  \cdots   \\[1.5pt]
				   +\: (1+1+1+1+\cdots+1_{\scriptscriptstyle\lambda})_{\scriptscriptstyle\lambda}
  \end{pmatrix}
   dx^2 \\[1.5pt]
  \end{split}
\end{align}
where evidently {\small$dx=dy=b/\La$}.

But now, just for a brief moment, let's forget about {\small$\La$}. For the very last time, we'll use the common {\small$+\cdots$} notation to denote an infinite ending, similarly to 1.\ref{eq:1} and 1.\ref{eq:2}. Therefore, if we add the terms in diagonals, starting from the top left corner, the double integral becomes 
\begin{align}  \label{eq:b2}
\begin{split}  
\int_0^b \int_0^b dxdy &= \Big( 1 + (1+1)+(1+1+1) +(1+1+1+1)+ \cdots \Big) \cdot dx^2  \\[2pt] 
				   &= \big( 1+2+3+4+ \cdots \big) \cdot dx^2 \\[3pt]
\end{split} 
\end{align}
However, if we now compare the Riemann sums in 1.\ref{eq:1} and 2.\ref{eq:b2}, the problematic conclusion would be that they're the same \textbf{---} and must thus be equal. But that's simply not true; the integral in 1.\ref{eq:1} gives {\small$b^2/2$}, while we obtain {\small$b^2$} in 2.\ref{eq:b2}. More specifically,  \vspace{0.6mm}
\begin{equation}
\int_0^b \int_0^bdxdy = 2\int_0^b xdx \\[4.5pt]
\end{equation} 
To resolve this would-be paradox, there's only one explanation: these Riemann sums have the same beginning, but different endings. In fact, by using {\small$\La$} again, we can demonstrate that their lengths are not the same; the Riemann sum {\small$(1+2+3+\cdots +\La) \cdot dx^2 $} contains {\small$\La$} terms, while there are {\small$2\La$} terms in 2.\ref{eq:b2}.

But, to understand this, we must first count the number of terms in the following double sum
\begin{equation} \label{eq:bbb3}
\sum^\La_{m=1} \sum^\La_{n=1} f(m,n) \\[4.5pt]
\end{equation}
where {\small$f(m ,n)$} can be any function. To do this, however, we need to extend the set {\small$\Na_{2\La}$} even further.  We thus define {\small$2\La+1$}, then {\small$2\La+2$}, and so on, up to the infinite limit {\small$2\La+\La=3\La$}. Similarly, we can continue with {\small$4\La$}, and {\small$5\La$}, and so on, up to yet another infinite limit; this gives {\small$\La\cdot\La=\La^2$}. Although we can keep counting even more, there's already a one-to-one correspondence between the resulting set {\small$\Na_{\La^2}$} and the terms in 2.\ref{eq:bbb3}. In particular, the sum with the index {\small$m$} contains {\small$\La$} terms, and each of them is yet another sum of {\small$\La$} terms \textbf{---} this gives a total of {\small$\La\cdot\La = \La^2$} terms. 

Evidently, the Riemann sum in 2.\ref{eq:b1} has also the same number of terms, as it's even forming an infinite square of length {\small$\La$}. Its geometry can actually explain what happens when we're adding its terms in diagonals. Simply consider the figure below.
\begin{center}
\textbf{$4 \times 4$ :} \:  \includegraphics[width=2.8cm, height=2.2cm]{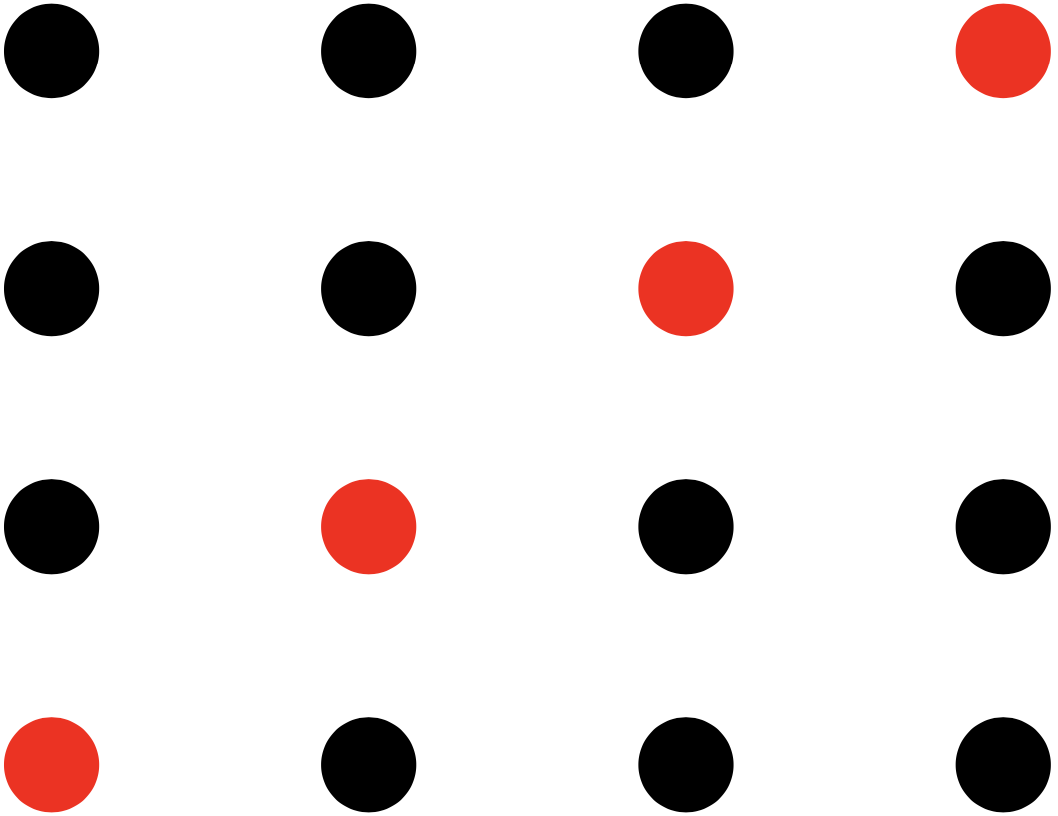} \quad\;\;\; \textbf{$\La \times \La$ :} \:  \includegraphics[width=2.8cm, height=2.2cm]{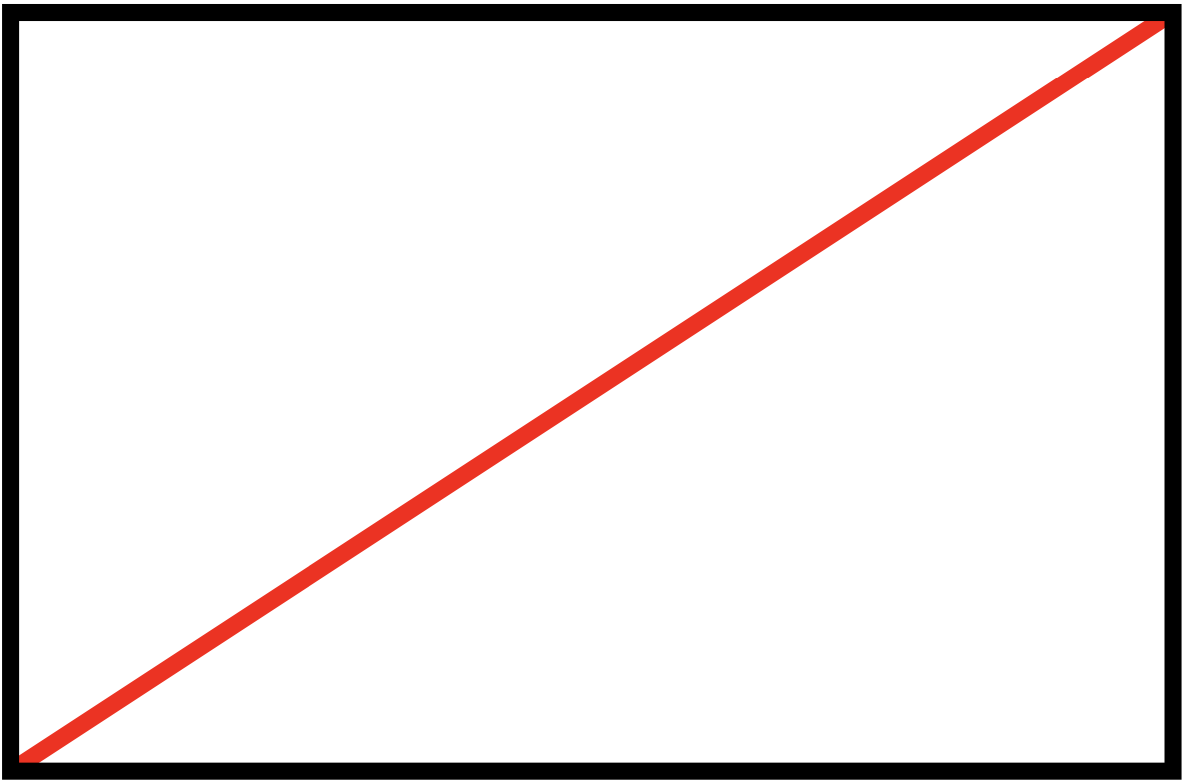} \quad\quad \\[9.5pt]
\end{center} 
In the finite square, the diagonals have successively $1$ dot, $2$ dots, then $3$, and the red one has $4$ dots; the remaining diagonals go down with $3$, $2$ and $1$ dots. Similarly, in the square of length {\small$\La$}, we start counting the number of dots in each diagonal with $1$, $2$, $3$, and so on. However, in this case, we must take an infinite limit to reach the red line, which now contains {\small$\La$} dots. There's thus a direct one-to-one correspondence between {\small$\Na$} and these diagonals: the $n$th one contains $n$ dots. To enumerate all the remaining diagonals, we must therefore use another copy of {\small$\Na$}. Starting at the bottom right corner, there are $1$, $2$, $3$ dots, and so on, where a second infinite limit must be taken to complete the enumeration. 

In the next section, we'll start by evaluating the number of dots in these infinite triangles. Similarly to the finite $4$x$4$ square, the second half will have one less diagonal than the first \textbf{---} we'll see this from the operation of subtraction. Nevertheless, for our current argument, we can safely assume that these triangles are the same. In fact, if we now replace all the dots in the {\small$\La$}x{\small$\La$} square with $1$s, and multiply them by {\small$dx^2$}, we'll recover the double Riemann sum in 2.\ref{eq:b1}. Finally, by adding those $1$s in diagonals, we get two sums {\small$(1+2+3+\cdots+\La)\cdot dx^2$}. 

Now, if we compare again the integrals in 2.\ref{eq:0}, but this time with our own {\small$\La$}-notation, we correctly obtain
\begin{align}
\begin{split}  \label{eq:bsym}
\int_0^b \int_0^b dxdy &= \: \big( 1+2+3+\cdots+\La + \La + \cdots + 3+2+1 \big) \cdot dx^2 \\[5pt]
				&= \: 2 \cdot \big( 1+2+3\cdots+\La \big) \cdot dx^2 \\[5pt]
				&= 2 \cdot \int_0^b x dx \\[0.5pt]
\end{split}
\end{align} 
The first conclusion is thus that {\small$\La<\La^2$}. But since we can write {\small$\La^2=\La\cdot\La$}, and since {\small$n<\La$} for any finite number $n$, the general inequality should rather be 
\begin{equation} 
\La \leq n\La < \La^2 \\[2pt]
\end{equation}

Needless to say, a valid theory should be able to provide plenty of such evidences. For this very purpose, we'll now turn to the widely used geometric series
\begin{equation} \label{eq:b3}
 1+x+x^2+x^3 + \cdots + x^n   = \frac{x^{n+1} - 1}{x-1}  \\[3.5pt]
\end{equation}
Typically, when {\small$n=\La$} and {\small$| \, x \, | < 1$}, the term {\small$x^{\La+1}$} is just assumed to be zero. The geometric series is thus commonly written as \vspace{0.5mm}
\begin{equation} \label{eq:badgeo}
 1+x+x^2+x^3 + \cdots = \frac{1}{1-x} \\[2pt]
\end{equation}
But that's a dangerous approximation: when {\small$x$} is infinitely close to {\small$1$}, the term {\small$x^{\La+1}$} is quite far from being zero. In particular, if we choose {\small$x = 1-dx$} with the infinitesimal {\small$dx=b/\La$}, and if we use the well-known Euler identity \vspace{1.8mm}
\begin{align}
\begin{split} \label{eq:euler}
(1-dx)^a = e^{a\ln(1-dx)} \sim e^{-adx}, \\[5pt]
\end{split}
\end{align}
then the geometric series in 2.\ref{eq:b3} becomes \vspace{0.8mm}
 \begin{equation} 
1 + e^{-dx} + e^{-2dx} + e^{-3dx} + \cdots + e^{-\La dx}  =  \frac{e^{-b} - 1}{-dx}  \\[4pt]
 \end{equation}
Now, by multiplying both sides with $dx$, the LHS becomes a Riemann sum, and we finally get \vspace{0.4mm}
 \begin{equation}
 \int_0^b e^{-x} dx = 1 - e^{-b}  \\[6pt]
 \end{equation}
 The term {\small$x^{\La+1}$} thus contributes to {\small$e^{-b}$} in the integral's value, and removing it from 2.\ref{eq:b3} would inevitably break this equality \textbf{---} even if {\small$| \, x \, | < 1$} is respected. 

The same problem arises, in fact, when people commonly write \vspace{0.8mm}
\begin{equation} \label{eq:wronggeo}
1 + 2x + 3x^2 + 4x^3 + \cdots = \frac{1}{(1-x)^2} \\[2pt]
\end{equation}
after squaring both sides of 2.\ref{eq:badgeo}. In particular, if {\small$x=1-dx$} is chosen again, this equality would break too. Multiplying both sides by {\small$dx^2$}, we'd get
\begin{align}
\begin{split}
\sum_{n=1}^\La n(1-dx)^{n-1}dx^2 &= \sum_{n=1}^\La  \: (ndx)  e^{-(n-1)dx}   dx \\[5pt]
						&\sim \int_0^b xe^{-x}dx \: \neq \: 1  \\[3pt]
\end{split}
\end{align}
where {\small$\sim$} is simply to get rid of the factor {\small$e^{dx}$}, since it's almost equal to $1$. But, fortunately, this problem is easy to fix. 

First, we must square both sides of 2.\ref{eq:b3} and expand {\small$\big(1+x+x^2 + \cdots + x^\La \big) ^2$} as a {\small$\La$}x{\small$\La$} square. Then, we add the terms in diagonals. Similarly to the double Riemann sum in 2.\ref{eq:bsym}, this gives us two infinite sums; that is, \vspace{0.45mm}
\begin{equation} \label{eq:rightgeo}
 \sum_{n=1}^{\La+1}  n x^{n-1} + \sum_{n=1}^\La \: n x^{2\La +1 - n} \: =  \bigg ( \frac{x^{\La+1} - 1}{x-1} \bigg ) ^ 2 \\[4.5pt]
\end{equation} 
The problem is now clear: there's a missing infinite sum in 2.\ref{eq:wronggeo} \textbf{---} that's the right half of our {\small$\La$}x{\small$\La$} square. By thus choosing {\small$x = 1-dx$}, and by multiplying both sides with {\small$dx^2$}, we can finally get the valid result \vspace{0.25mm}
\begin{align}
\begin{split}
\quad\;\;  \sum_{n=1}^{\La+1} n(1 - dx)^{n-1} dx^2 + \sum_{n=1}^\La n (1 - dx)^{2\La + 1 - n} dx^2 \quad\quad\quad\quad\quad \\[6pt]
\sim \int_0^b xe^{-x}dx + e^{-2b} \int_0^b xe^x dx \quad\quad\quad\quad \\[11pt]
= (e^{-b} - 1)^2 \quad\quad\quad\quad\quad\quad\quad\quad\quad\quad\quad\quad  \\[1.5pt]
\end{split}
\end{align}

The identity 2.\ref{eq:wronggeo} is thus demonstrably incomplete; it's missing the {\small$\La$} terms that gave us {\small$ e^{-2b} \int_0^b xe^x dx$}. That's yet further evidence for the existence of {\small$\La^2 =\La\cdot\La$} and its square-like structure,  which is notably much bigger than {\small$|\Na|$}. At this point, it should be fairly intuitive that {\small$\La^n$} must exist too, with the inequality 
\begin{equation}
\La^m<\La^n \quad \text{for} \quad m<n \\[2pt]
\end{equation}
The upcoming sections will use such quantities extensively, and we'll thus gather many more evidences for their existence.

Moreover, we also saw that infinitesimal quantities can have a major impact in calculations \textbf{---} that is, when infinitely many of them are added together. Discarding them can often bring welcomed simplifications, but it must be done with great caution. Adding {\small$\La$} of an infinitesimal quantities can lead to different results, depending on its size. For a quick example, \vspace{0.2mm}
\begin{equation}
\La\cdot(1/\La^2) = 1/\La\sim0 \quad\quad \text{but} \quad\;\;\: \La\cdot(1/\La)=1 \quad\:
\end{equation}

Now, before ending this section, we must evaluate the cardinality of {\small$\QQ$}, the set of all rational numbers. That's easy to achieve, since they're written as the ratio {\small$a/b$}, where {\small$a,b \in \Na$}. There are thus {\small$\La$} possible numerators, and, for each of them, there are as many denominators. We therefore already have
\begin{equation}
\;\;\: |\QQ| = \La \cdot \La = \La^2 \\[0.1pt]
\end{equation}

Unsurprisingly, our result is yet again in direct opposition to Cantor's theory, where the sets {\small$\Na$} and {\small$\QQ$} are said to have the same cardinality. This would result from the alleged one-to-one correspondence \vspace{0.3mm}
\begin{center}
\quad\quad\quad \includegraphics[width=4.6cm, height=3.65cm]{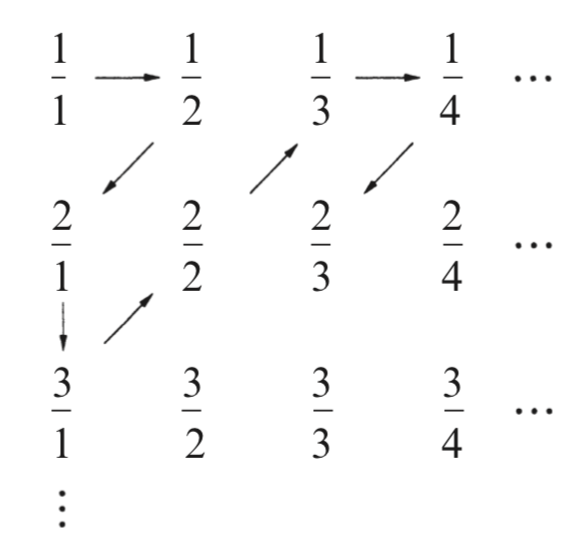}
\end{center}
Nevertheless, similarly to Cantor's mapping between {\small$\Na$} and {\small$\ZZ$}, this is completely wrong \textbf{---} and, in this case, it's even easier to see why. 

In the figure above, the rational numbers are clearly forming an infinite square of length {\small$\La$}, and Cantor's mapping is zigzagging through its diagonals. However, even if we map the $n$th number of {\small$\Na$} to the $n$th diagonal, with all the rational numbers that it contains, this wouldn't be a one-to-one correspondence. There are {\small$2\La$} diagonals in this square, and, to enumerate all the diagonals, we must use another copy of {\small$\Na$}. 

To give a concrete example, the rational number {\small$\frac{(\La-3)}{(\La-2)}$} is left out of Cantor's mapping. But it's a meaningful member of {\small$\QQ$}, because we can multiply it by itself {\small$\La$} times, and, from 2.\ref{eq:euler}, the result will be \vspace{0.8mm}
\begin{equation} \label{eq:rational}
\Bigg( \frac{\La-3}{\La-2} \Bigg)^\La = \Bigg( \frac{1-3/\La}{1-2/\La} \Bigg)^\La =  \frac{\big(1-3/\La\big)^\La}{\big(1-2/\La\big)^\La} = \frac{e^{-3}}{e^{-2}} = e^{-1} \\[4pt]
\end{equation}

Admittedly, we cheated a little by using both subtraction and exponentiation in 2.\ref{eq:rational}, since they will be introduced, respectively, in Sections $3$ and $5$. Nevertheless, the ordering {\small$\La \leq n\La < \La^2$} is now firmly established in our work, and this continues to challenge Cantor's theory \textbf{---} as it was shown, once again, to be inconsistent with Calculus. \\[9pt]

 \setcounter{equation}{0}
 \setcounter{section}{3}
 \vspace{1mm}

\section*{3. Subtraction }
So far, {\small$\La$} was used in additions (i.e. {\small$\La+n$}) and multiplications (i.e. {\small$n\La$} and {\small$\La^n$}). It's therefore time to introduce an inverse operation \textbf{---} namely, subtraction. Similarly to how we start from $1$ and take an infinite limit to reach {\small$\La$}, we can also start from {\small$\La$} and make our way down to $1$; we simply define {\small$\La-n$}. 

Surprisingly perhaps, our first use of subtraction will be to evaluate infinite sums of numbers. Let's begin with the simplest one: the sum of all {\small$\Na$}'s numbers. That is, we want to evaluate 
\begin{equation}
S_1 = \sum^\La_{n=1} n  \\[3pt]
\end{equation}
Graphically, the sum is 
\begin{center}
\quad\quad \: \includegraphics[width=4.5cm, height=2.8cm]{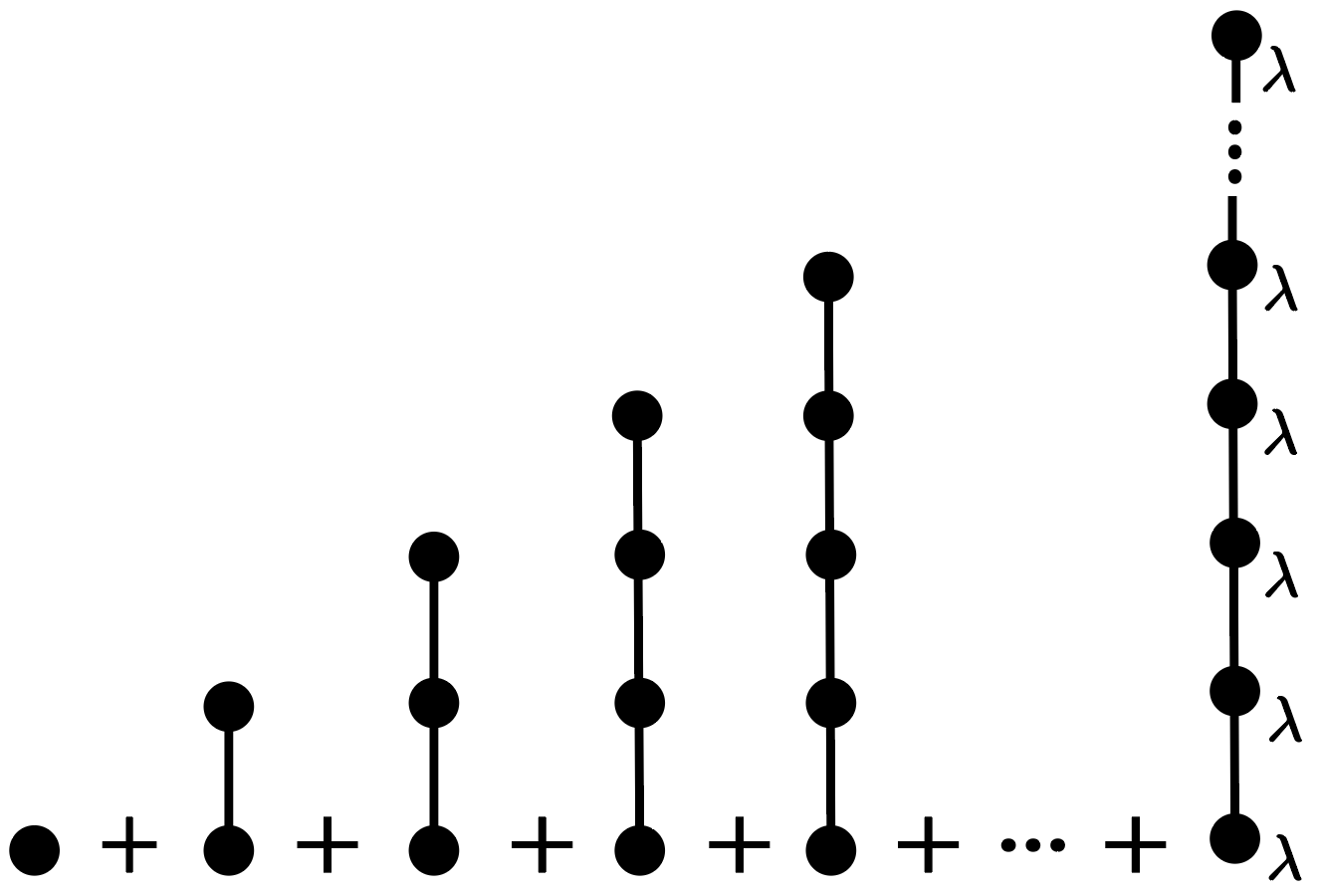} \quad\quad\quad\quad\quad\quad 
\end{center} \vspace{0.4mm}
where each {\small$\Na$}'s numbers is shown as a column of black dots. Clearly, to evaluate {\small$S_1$}, we must simply count all these dots. That's easier to achieve if we group them horizontally instead, since it leads to the new visualization \vspace{1.5mm}
\begin{center}
\quad\;\:\: \includegraphics[width=4.6cm, height=3cm]{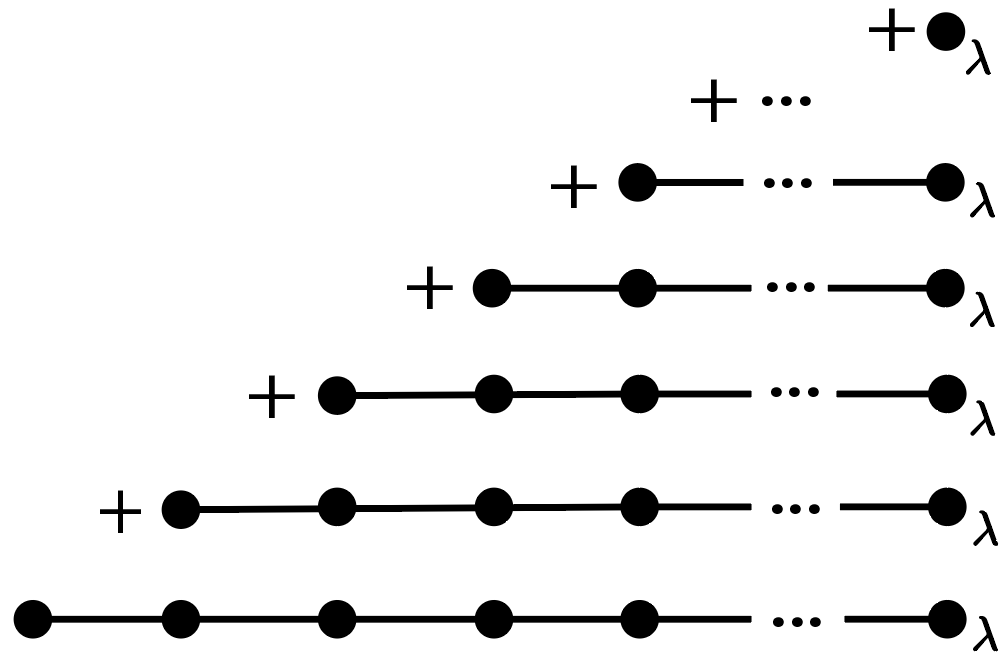} \quad\quad\quad\quad\quad\quad
\end{center} \vspace{0.4mm}
At the very bottom, there's a line with {\small$\La$} dots, one for each number in {\small$\Na$}. In the next line, however, there's one less dot; we thus count {\small$\La-1$} of them. Going through all the other lines, we have successively {\small$\La-2$} more dots, then {\small$\La-3$}, and so on. Adding them all together, we get the total count
\begin{align}
\begin{split}
S_1 &= \sum_{n=1}^{\La} \big(\La - n +1\big) \\[5.5pt]
 &   = \La^2 - S_1 + \La  \\[2pt]
\end{split}
\end{align}
Solving for {\small$S_1$}, we finally get 
\begin{equation} \label{eq:c1}
S_1 = \frac{\La(\La+1)}{2} \\[3.5pt]
\end{equation}
This formula should be familiar, as it's also valid when {\small$S_1$} is finite. Moreover, we can now verify the number of dots in the {\small$\La$}x{\small$\La$} square. In particular, since the left triangle has {\small$\La$} diagonals, and since the right one has {\small$\La-1$}, the total number of dots must be \vspace{0.2mm} 
\begin{equation}
\sum_{n=1}^{\La} n +  \sum_{n=1}^{\La-1} n \: = \: \frac{\La(\La+1)}{2} + \frac{(\La-1)(\La)}{2} \: = \: \La^2 \\[5pt]
\end{equation}

Now, for yet another infinite sum, we'll add the first {\small$\La$} triangular numbers. That is,
\begin{equation}
S_2 = \sum_{n=1}^\La \frac{(n)(n+1)}{2} \\[4pt]
\end{equation} 
where the sum starts with the terms $1$, $3$, $6$, $10$, $15$, and so on; the {\small$\La$}th one is actually {\small$S_1 = \frac{\La(\La+1)}{2}$}. The procedure will be the same: we first count {\small$\La$} by removing $1$ from each term in the above sum, which gives us 
\begin{equation}
S_2 = \La + \Big(0+2+5+9+14+\cdots+(S_1-1)\Big) \\[2pt]
\end{equation}  
Since the first non-zero term is now {\small$2$}, we can count another {\small$2(\La-1)$} by removing {\small$2$} to all the remaining terms. This leads to 
\begin{equation}
S_2 = \La + 2(\La-1) + \Big(0+0+3+7+12+\cdots+(S_1-3) \Big)
\end{equation}
Next, we count another {\small$3(\La-2)$ by removing $3$ to all the non-zero terms. By doing this until the initial sum is empty, we ultimately get   
\begin{align}
\begin{split} \label{eq:33}
S_2 &= \sum_{n=1}^{\La} n\big(\La-n+1\big)  \\
    &= \sum_{n=1}^{\La} n\big(\La + 2 - (n+1) \big) \\[4.5pt]
    &=S_1(\La+2) - 2S_2 \\[5pt]
\end{split}
\end{align} 
Solving for {\small$S_2$}, and using 3.\ref{eq:c1}, the final result is \vspace{0.6mm}
\begin{equation} \label{eq:d2}
S_2 = \frac{(\La)(\La+1)(\La+2)}{3!}  \\[2pt]
\end{equation} 
Once again, this is valid when {\small$S_2$} is finite. By using induction with the same procedure, it's now straightforward to generalize to the familiar  \vspace{0.9mm}
\begin{equation}
S_b = \sum_{n=1}^{\La} \binom{n+b}{n} = \frac{(\La)(\La+1)(\La+2)\cdots(\La+b-1)}{b!} \\[2pt]
\end{equation}

Across all of science, an infinite sum of numbers is always equal to the same {\small$\infty$}, and any useful information is thus destroyed irreversibly. But it doesn't have to be this way. An infinite sum simply keeps the same closed-form expression as when it's finite; there's no reason why taking an infinite limit would erase any of its content. In fact, we'll shortly make great use of this \textbf{---} although, for the moment, that's not even the main takeaway. Similarly to {\small$\La+n\neq\La$}, we now also have the inequality {\small$\La-n\neq\La$}. For example, by writing {\small$\La-2$}, we refer to the sum {\small$0+0+1+1+\cdots+1_\sLa$}, which is not equal to {\small$\La=1+1+1+\cdots+1_\sLa$}. 

As before, we'll use Calculus to verify that our results are valid. For this purpose, consider the following integral and its Riemann sum.
\begin{align}
\begin{split} \label{eq:832}
\int_0^b  \big(b-x \big)xdx &= \sum_{n=1}^\La \big(b-ndx \big)(n dx)dx \\[3.25pt]
			      &= \sum_{n=1}^\La   \big(\La-n \big) n dx^3 \\[1pt]
\end{split} 
\end{align}
To simplify the notation a little bit, we'll use the Riemann sum {\small$\sum^\La_{n=1} (\La-n+1)ndx^3$} instead \textbf{---} it's equal to 3.\ref{eq:832} plus only one infinitesimal term. Now, by using our interpretation for {\small$\La-n$}, we can continue with \vspace{0.2mm}
\begin{equation}
\int_0^b (b-x)xdx =  \begin{pmatrix}
			
			\:\:\:\:\:1\cdot(\La-0) \\[1pt]
			+\: 2\cdot(\La-1) \\[1.5pt]
			+ \:3\cdot(\La - 2) \\[1.5pt]
			+\: 4\cdot(\La - 3) \\[2pt]
			 \cdots \\[1.5pt]
			+\: \La \cdot (1) \:\:\:\:\:\:\:\:\:

			\end{pmatrix} dx^3
=  		 \begin{pmatrix}
 			\:\:	\:\:\:1+1+1+1+\cdots+1_\sLa \\[1.5pt]
				   + \: 0+2+2+2+\cdots+2_\sLa \\[1.5pt]
				   + \:0+0+3+3+\cdots+3_\sLa \\[1.5pt]
				+ \:0+0+0+4+\cdots+4_\sLa \\[2pt]
				   \;\:  \cdots   \\[1pt]
				  +\:  0+0+0+0+\cdots+\La_\sLa
  \end{pmatrix} dx^3 \\[3pt]
\end{equation}
Adding all the terms in a vertical line together, we successively have {\small$1$}, {\small$3$}, {\small$6$}, {\small$10$}, {\small$15$}, and so on. Finally, using 3.\ref{eq:d2}, we conclude with
\begin{align}
\begin{split} 
\int_0^b (b-x)xdx  &=   \sum_{n=1}^\La \frac{(n)(n+1)}{2} \cdot dx^3  \\[4pt]
			       &= \bigg( \frac{\La^3+3\La^2+2\La}{3!} \bigg) \cdot \bigg( \frac{b}{\La} \bigg)^3 \\[7pt]
			       &\sim b^3/3
\end{split}
\end{align} 
as expected. Clearly, to evaluate this Riemann sum correctly, our result {\small$\La\neq\La-n$} is required to be true.

For yet another way to validate this, we can resolve paradoxes \textbf{---} and, unsurprisingly, there are several that emerge from the widespread identity \vspace{0.2mm}
\begin{equation}
\sum_{x=1}^{\La} f(x) = \sum_{x=1}^{\La-n} f(x) \\[2pt]
\end{equation}
Typically, it leads to problematic results similar to the one below, where the geometric series is erroneously derive from
\begin{align}
\begin{split} 
s &= 1 + b + b^2 + b^3 + b^4 + \cdots \\
   &= 1 + b ( 1 + b^2 + b^3 + \cdots) \\
   &= 1 + bs
\end{split} 
\end{align}
By solving for $s$, we would obtain the common {\small$s=1/(1-b)$}. But if  {\small$b=2$}, we'll get the nonsensical identity {\small$s=-1$}, which is taken seriously by some mathematicians\cite{Gardiner}. Nevertheless, by using {\small$\La$} instead of the incomplete {\small$+\cdots$} notation, we get
\begin{align}
\begin{split} 
s & = 1 + b + b^2 + b^3  + \cdots + b^{\lambda} \\
   & = 1 + b(  1 + b + b^2  + \cdots +  b^{\lambda-1}) \\
   & = 1 + b(s - b^{\lambda})
\end{split} 
\end{align}
Solving for $s$, we obtain 2.\ref{eq:b3} back, with {\small$n=1$}. That's expected, since this is also valid for the finite case. 

The same problem arises with the famous claim that $1 = 0.999\ldots$, which was advocated by Leonhard Euler himself in \textit{Elements of Algebra} (1770). By following Euler's argument, but using the correct manipulations instead, we can find that
\begin{align}
\begin{split} 
\quad\quad\quad\quad\;\:  x &= 0.999\ldots \\
	&= 9/10 + 9/10^2 + 9/10^3 + \cdots + 9/10^\lambda  \\[4pt]
10x &= 9 + \big( 9/10 + 9/10^2 + \cdots + 9/10^{\lambda-1} \big)  \\[4pt]
10x &= 9 + \big( x - 9/10^\lambda \big) \\[4pt]
  x & = 1 - 1/10^\lambda
\end{split} 
\end{align}
Mathematicians would typically argue that infinitesimal quantities, such as {\small$1/10^\La$}, are equal to zero. But that's simply not true: by adding infinitely many of them, we can compute an integral from its Riemann sum. 

In fact, if we multiply {\small$0.999\ldots$} by itself {\small$10^\La$} times, by using ironically the Euler identity 2.\ref{eq:euler}, we get {\small$(1-1/10^\La)^{10^\La} = e^{-1}$}. On the other hand, we can multiply $1$ by itself any number of times, and we'll always get $1$. Consequently, the assertion that {\small$1=0.999\ldots$} just can't be true. (Again, exponentiation with {\small$\La$} will be covered in Section $5$.)

There's yet an entire class of paradoxes that {\small$\La \neq \La-n$} can resolve. Generally, they're known as \textit{paradoxical decompositions}: an object is first decomposed into distinct pieces, and it's then reassembled into more copies of itself. For example, the Banach-Tarski paradox states that a solid sphere can be cut into finitely many pieces, and rearranged to give back two identical ones \cite{Banach-Tarski, wagon}. Nevertheless, since all these paradoxes have a very similar resolution, we'll focus on just one of them  \textbf{---}  chosen for its simplicity. 

In his book \textit{From Here to Infinity} (1996), Ian Stewart introduces the Hyperwebster, a dictionary that contains all possible words made from our 26-letter alphabet. The Hyperwebster is divided into volumes: all words beginning with A are in Volume A, those with B are in Volume B, and so on.  \vspace{2.4mm}

\quad \includegraphics[width=11cm, height=3.9cm]{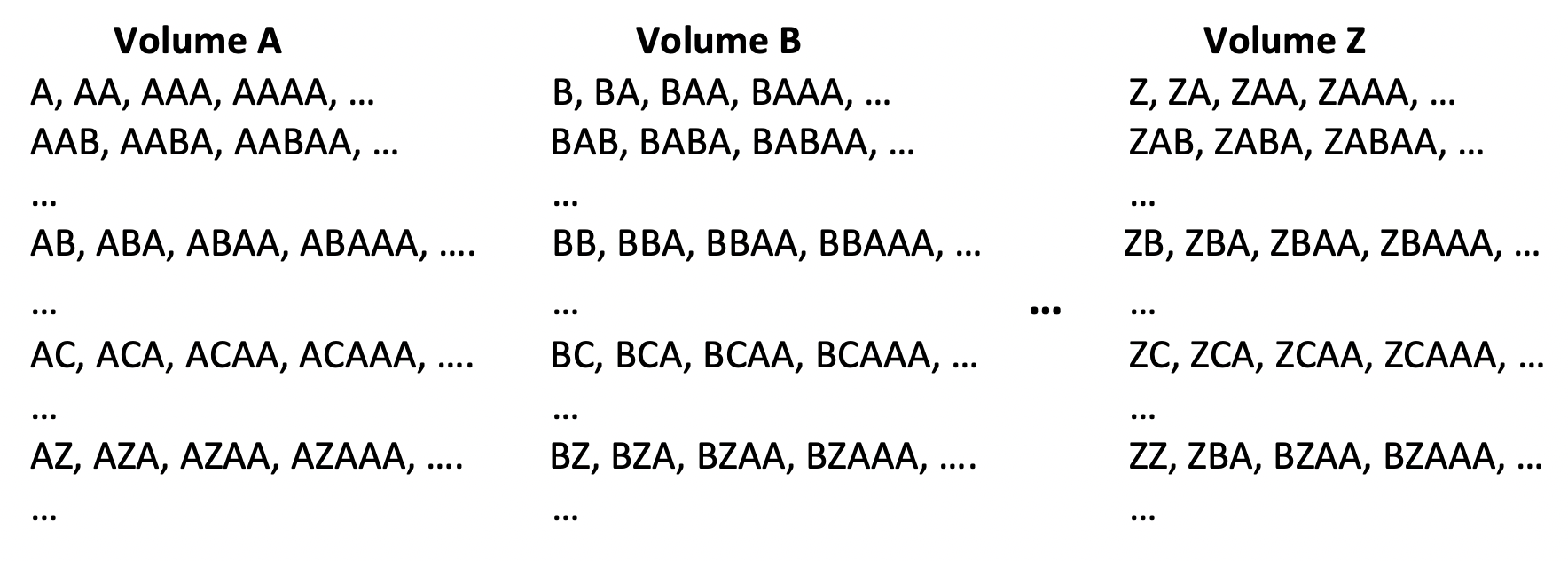}

According to Stewart, by removing the first letter in every word, we can produce another copy of the Hyperwebster from each of its volumes; that is, {\small$26$} identical copies will be generated. Doing this again on every resulting copy, we can now have {\small$26^2$} Hyperwebsters, then {\small$26^3$}, and so on, with no limits whatsoever. Amusingly, Stewart writes that "in spirit, the Banach-Tarski paradox is just the Hyperwebster wrapped round a sphere". 

Nevertheless, these paradoxical duplications, where objects are created out of nothing, are only made possible from a faulty believe \textbf{---} that removing finite quantities from infinity will leave it unchanged. However, by using {\small$\La \neq \La-n$}, we can resolve this paradox. In particular, all words with $n$ letters are erased after the $n$th duplication; it's thus clear that each Hyperwebster's copy is shorter than the original one. In fact, after the {\small$\La$}th duplication, if each word is at most {\small$\La$} letters long, we'll be left with {\small$26^{\La}$} empty dictionaries.

Finally, as the last paradox of this section, we must address the infamous \vspace{0.2mm}
\begin{equation} \label{eq:ram}
\sum^\infty_{n=1} n = -1/12 
\end{equation}
which is generally attributed to Srinivasa Ramanujan \cite{ramanujan}. Evidently, the correct identity is already given in 3.\ref{eq:c1}, and we could stop there. But this paradox has spread widely into the scientific community; even physicists are trying to use it in their work \cite{stringtheory, StandardModel, Bernardo}. Consequently, we want  to clearly explain what's wrong in Ramanujan's calculations, and hopefully this will prevent it from causing anymore damage. Before doing so, however, we have to present Euler's own ideas on the subject, since 3.\ref{eq:ram} heavily relies on them. 

In \textit{De seriebus divergentibus} (1760), Euler's goal was to evaluate a few infinite sums, with {\small$1-2+3-4+\cdots$} in particular. But he had to first recognize a serious shortcoming: since the terms are alternatively positive and negative, the partial sums also follow the same pattern; it's a divergent series. To circumvent this issue, Euler had the idea to use the geometric series, and he wrote \vspace{0.2mm}
\begin{equation}
1-2x+3x^2-4x^3+\cdots = \frac{1}{(1+x)^2} \\[2.5pt]
\end{equation}
Recognizably, this is already an error, because the valid identity is 2.\ref{eq:rightgeo}. Euler also stated that {\small$| \, x \, | < 1$} is a requirement for this to hold, but went on to simply pick {\small$x=1$}. Most mathematicians now agrees with his method, although they insist on using the limit of $x$ approaches $1$, and call it Abel's theorem \cite{hardy}. In fact, they conclude with 
\begin{equation} \label{eq:eulerparadox}
\lim_{x\rightarrow 1^{-}} \sum_{n=1}^\infty n(-x)^{n-1} = \lim_{x\rightarrow 1^{-}}  \frac{1}{(1+x)^2} = 1/4 
\end{equation}
Something must be wrong, however; adding and subtracting integers can't possibly lead to a fraction. Before continuing any further, we thus have to evaluate 3.\ref{eq:eulerparadox} correctly.

This series is divergent for us too, however, and we can't evaluate it with only {\small$\La$} terms \textbf{---} we'd be forced to use {\small$\lim_{n\rightarrow\La} (-1)^n = (-1)^\La$}, which is also divergent. To avoid this problem, we'll simply take this sum with instead {\small$2n\La$} terms, where $n \in \Na$. We can thus use 
\begin{equation}
(-1)^{2n\La} = \Big( (-1)^{2n} \Big)^\La = 1^\La = 1 \quad \text{and} \quad (-1)^{2n\La+1} = -1
\end{equation}
In fact, we can already conclude with \vspace{0.4mm}
\begin{align} \label{eq:rightram}
\begin{split}
1-2+3-4+\cdots-2n\La &= \Big(1+3+5+\cdots+(2n\La-1) \Big) \quad\quad\quad\quad\;\: \\[0.2pt]
&\quad\quad\; - \Big(2+4+6+\cdots+2n\La\Big) \\[3pt]
						&= \sum_{j=1}^{n\La} (2j-1) -  \sum_{j=1}^{n\La} (2j) \\[2.5pt]
						&= \sum_{j=1}^{n\La} (-1) = -n\La \\[1pt]
\end{split} 
\end{align}

Now, to finally resolve Ramanujan's paradoxical identity, we must first explain how he got there. He started by writing 
\begin{alignat}{7}  \label{eq:wrongram}
 c&{}={}&1+2&&{}+3+4&&{}+5+6+\cdots \\
4c&{}={}&  4&&  {}+8&&{} +12+\cdots  \label{eq:wrongram} \\
c-4c&{}={}&1-2&&{}+3-4&&{}+5-6+\cdots \label{eq:wrongram2} 
\end{alignat}
From Euler's own paradoxical identity 3.\ref{eq:eulerparadox}, Ramanujan then claimed that \vspace{0.55mm}
\begin{equation}
-3c = 1-2+3-4+5-6+\cdots = 1/4 \\[0.2pt]
\end{equation}
The final result would thus be {\small{$c=-1/12$}. Nowadays, mathematicians insist on doing these manipulations with the zeta function and a limit of $x$ approaching $1$, similarly to 3.\ref{eq:eulerparadox}. But, fundamentally, this is just as wrong.

To resolve this paradox, and due to the divergence problem, we'll take $c$ to have {\small$2\La$} terms. Now, in 3.\ref{eq:wrongram}, notice that Ramanujan is actually writing \vspace{0.2mm}
\begin{equation}
4c=0+4+0+8+0+12+\cdots  \;\;
\end{equation}
This sum therefore contains {\small$4\La$} terms, because there are {\small$2\La$} zeroes and {\small$2\La$} numbers. But that's a serious issue in 3.\ref{eq:wrongram2}: the last {\small$2\La$} terms in $4c$ do not subtract any term of $c$, due to their difference in length. Since there are {\small$\La$} zeroes in these terms, we conclude that {\small$\La$} numbers are missing in Ramanujan's calculations. These are  \vspace{0.4mm}
\begin{equation}
4 \cdot \Big( c - (1+2+3+\cdots+\La ) \Big) \\[1pt]
\end{equation} 
Ramanujan thus made two mistakes: first, using Euler's identity 3.\ref{eq:eulerparadox} is wrong; and, due to the added zeroes, the subtraction in 3.\ref{eq:wrongram2} is not one-to-one. 

Now, to verify this, we just need to evaluate $c-4c$ with our own results, and solve for $c$ to see if we get a valid outcome. From 3.\ref{eq:c1} and 3.\ref{eq:rightram} (with {\small$n=1$}), we thus calculate 
\begin{align}
\begin{split}
-3c &= \big(1-2+3-4+ \cdots - 2\La \big) - 4 \cdot \bigg(c - \sum_{n=1}^\La n \bigg) \\[3pt]
	&= -\La - \Big( 4c - 2 \big( \La^2 + \La \big) \Big) \\[4.5pt]
\end{split} 
\end{align}
\pagebreak
and this finally leads to {\small$c=2\La^2+\La$}. That's clearly the expected result, since we could get it by simply setting {\small$b=2\La$} in the well-known formula \vspace{0.4mm}
\begin{equation}
\sum^b_{n=1} n = \frac{b(b+1)}{2} \\[1pt]
\end{equation}

It's hard to grasp why so many mathematicians prefer to accept nonsensical results such as Ramanujan's, rather than step back and look for errors.
 By using the true properties of infinity, however, we could find that infinitely many terms were missing in Ramanujan's calculations \textbf{---} and this lead to the unfortunate identity 3.\ref{eq:ram}. In fact, in the next section, we'll demonstrate that even the great Bernhard Riemann did a very similar mistake. \\[12pt]

\setcounter{equation}{0} 
\setcounter{section}{4}

\section*{4. Division} 
\vspace{0.6mm}
After studying subtractions, we can now explain how it's possible to divide {\small$\La$}. Specifically, we want to use and interpret expressions such as {\small$\La/n$}. For this very purpose, we start with \vspace{0.6mm}
\begin{equation}
\La = 1+1+1+\cdots+1_{\scriptscriptstyle\lambda} \\[3.7pt]
\end{equation}
from which it's natural to define \vspace{1.4mm}
\begin{equation}
\frac{\La}{n} = \frac{1}{n}+\frac{1}{n}+\frac{1}{n}+\cdots+\frac{1_{\scriptscriptstyle\lambda}}{n} \\[2.8pt]
\end{equation}
If the first  $n$ terms are added together, they will sum up to $1$, and only one term will remain. Doing this again with the next $n$ terms, and so on, we can reduce our initial {\small$\La$} terms to the desired number of them \textbf{---} that is, {\small$\La/n$}.

To give a concrete example of such a quantity, we'll start by counting the number of odd numbers in the set {\small$\Na_n=\{1,2,3,...,n\}$}.  We either get exactly {\small$n/2$}, or very close to it, depending on how large $n$ is. By now taking the limit of $n$ goes to {\small$\La$}, we obtain that {\small$\Na$} contains {\small$\La/2$} odd numbers (or infinitely close to it). Since the same applies for the even numbers, it's safe to write that  {\small$|\Na| \sim \La/2+\La/2$}.

Admittedly, it's all very simple. But this now allows us to resolve a major paradox in mathematics: the so-called \textit{Riemann Rearrangement Theorem} \cite{Riemann, tao}. Loosely stated, it says that if an infinite series is convergent, then its terms can be rearranged so that the new series will converge to an arbitrary real number. A very famous example, also attributed to Peter Lejeune-Dirichlet (1827), goes as follows. First, the power series of {\small$\ln(2)$} is divided into three other ones.
\begin{align}
\ln(2) &= 1 - \frac{1}{2}+\frac{1}{3}-\frac{1}{4}+\frac{1}{5}-\cdots \label{eq:ln}  \\[4pt]
	&= 1+\frac{1}{3}+\frac{1}{5}+\frac{1}{7}+\frac{1}{9}+\cdots  \label{eq:e1}  \\[3pt]
	&-\frac{1}{2}-\frac{1}{6}-\frac{1}{10}-\frac{1}{14}-\frac{1}{18}-\cdots \label{eq:e2}  \\[4pt]
	&-\frac{1}{4}-\frac{1}{8}-\frac{1}{12}-\frac{1}{16}-\frac{1}{20}-\cdots \label{eq:e3}
\end{align}
 \\[0.1pt]The series 4.\ref{eq:e2} is then added to 4.\ref{eq:e1}, as to give
\begin{equation}
\bigg(1-\frac{1}{2}\bigg)+\bigg(\frac{1}{3}-\frac{1}{6}\bigg)+\bigg(\frac{1}{5}-\frac{1}{10}\bigg) + \bigg(\frac{1}{7}-\frac{1}{14}\bigg) + \cdots \label{eq:e4} \\[3pt]
= \frac{1}{2}+\frac{1}{6}+\frac{1}{10}+\frac{1}{14}+ \cdots  
\end{equation}
Finally, the series 4.\ref{eq:e3} is added to 4.\ref{eq:e4}, and this leads to the paradoxical conclusion 
\begin{align}
\begin{split}  \label{eq:e5}
\ln(2) &= \frac{1}{2}-\frac{1}{4}+\frac{1}{6}-\frac{1}{8}+\frac{1}{10}-\frac{1}{12}  + \cdots  \\[6pt]
          &= \frac{1}{2} \cdot \bigg(1 - \frac{1}{2}+\frac{1}{3}-\frac{1}{4}+\frac{1}{5}-\frac{1}{6}+\cdots\bigg) \\[6pt]
          &=  \frac{1}{2}\cdot\ln(2)
\end{split} 
\end{align}

Evidently, this must be wrong, since Riemann's theorem is saying that $2=1$. From our work, however, it's easy to see where's the error. Without loss of generality, we'll assume that 4.\ref{eq:ln} is made of {\small$\La$} terms \textbf{---} in other words, we'll write {\small$\ln(2) = \sum_{n=1}^\La (-1)^{n+1}/n$}. This now requires 4.\ref{eq:e1} to only have {\small$\La/2$} terms, since that's how many odd numbers there are in {\small$\Na$}. The other {\small$\La/2$} terms are equally divided between 4.\ref{eq:e2} and 4.\ref{eq:e3}, so they each have {\small$\La/4$} terms. 

Riemann therefore repeats twice the same mistake. When adding 4.\ref{eq:e2} to 4.\ref{eq:e1}, and when adding 4.\ref{eq:e3} to 4.\ref{eq:e4}, he doesn't realize that those series have different lengths; in both cases, there's one with {\small$\La/4$} terms while the other has {\small$\La/2$} terms. But since all terms in the longer series are expected to receive a term from the shorter one, there's a total of {\small$\La/4+\La/4=\La/2$} terms that Riemann did not include in his calculations. Therefore, if those missing terms are equal to {\small$\ln(2)/2$}, we can restore the equality in 4.\ref{eq:e5} \textbf{---} and the paradox will be resolved.

We start by evaluating the missing terms in 4.\ref{eq:e3}, a series that can now be written as {\small$\sum_{n=1}^{\La/4} \frac{1}{4n}$}. Because the last term is currently {\small$\frac{1}{4\cdot(\La/4)} = \frac{1}{\La}$}, for this series to be any longer, we need {\small$\frac{1}{\La+4n}$} to be the next terms. As a result, the {\small$\La/4$} missing terms to evaluate are  
\begin{equation}
\sum_{n=1}^{\La/4} \frac{1}{\La + 4n} =  \frac{1}{\La} \cdot \sum_{n=1}^{\La/4}  \frac{1}{1 + 4n/\La} \label{eq:e6} \\[5.5pt]
\end{equation} 

There are two ways to compute this: we can either use \textbf{---} ironically \textbf{---} Riemann's definition of an integral, or we can expand {\small$\frac{1}{1 + 4n/\La}$} as a geometric series. For the first method, we simply start by considering the integral
\begin{align}
\begin{split} 
\quad \quad \: \int_0^4 \frac{dx}{1+x} &= \sum^{\La}_{n=1} \frac{4/\La}{1 + 4n/\La}  \\[4.5pt]
				   &= \sum^{\La/4}_{n=1} \frac{4/\La}{1 + 4n/\La} + \sum^{\La}_{n=\La/4} \frac{4/\La}{1 + 4n/\La} \\[6.8pt]
				   &=  \int_{0}^1 \frac{dx}{1+x} + \int_{1}^4 \frac{dx}{1+x} \\[3pt]
\end{split}
\end{align}
In particular, if we compare 4.\ref{eq:e6} to the sum with {\small$\La/4$} terms, it's clear that they only differ by a factor of $4$. We can thus already conclude with
\begin{equation}
\frac{1}{\La} \cdot \sum_{n=1}^{\La/4}  \frac{1}{1 + 4n/\La} = \big(1/4\big) \cdot \int_{0}^{1}  \frac{dx}{1 + x} =  \frac{\ln(2)}{4} 
\end{equation} 

In the second method, each term {\small$\frac{1}{1 + 4n/\La}$} is expanded as a geometric series. This is a bit more messy, but we achieve the same result. That is,
\begin{align}
\begin{split} 
 \frac{1}{\La} \cdot \sum_{n=1}^{\La/4}  \frac{1}{1 + 4n/\La} \: &\sim \: \frac{1}{\La} \cdot \sum_{n=1}^{\La/4} \Bigg( \sum^\La_{b=0} \big(-4n/\La \big)^b  \Bigg) \\[5pt]
											&= \: \frac{1}{\La} \cdot \sum^\La_{b=0} \big( -4/\La \big)^b \sum_{n=1}^{\La/4} n^b \\[5pt]
											&= \:  \frac{1}{\La} \cdot \sum^\La_{b=0} \big( -4/\La \big)^b \cdot \Bigg( \frac{1}{b+1} \sum^b_{j=0} \binom{b}{j} B_j \big(\La/4\big)^{b-j} \Bigg) \\[5pt]
											&\sim\: \frac{1}{4} \cdot \sum^\La_{b=0} \frac{(-1)^b}{b+1} = \frac{\ln(2)}{4} \\[1pt]
\end{split} 
\end{align}
where the symbol $\sim$ was used to remove all the infinitesimal terms, and where {\small$B_j$} are the Bernoulli numbers. Obviously, by taking {\small$\La$} to be just a very big number, the above calculation can be verified numerically. 

Now, for the missing terms in 4.\ref{eq:e2}, we can quickly demonstrate that they give the same result. Since this series can be written as {\small$\sum_{n=1}^{\La/4} \frac{1}{4n-2}$}, the last term is evidently {\small$\frac{1}{4\cdot(\La/4)-2} = \frac{1}{\La-2}$}. Therefore, the missing terms are similarly given as
\begin{align}
\begin{split} 
\sum_{n=1}^{\La/4} \frac{1}{\La + 4n - 2} &=  \frac{1}{\La} \cdot \sum_{n=1}^{\La/4}  \frac{1}{1 + 4n/\La - 2/\La} \\[5.7pt]
							      &\sim \frac{1}{\La} \cdot \sum_{n=1}^{\La/4}  \frac{1}{1 + 4n/\La}  \\[7.5pt]
							      &= \frac{\ln(2)}{4} \\[2.4pt]
\end{split} 
\end{align} 
where the symbol $\sim$ was for removing the term {\small$-2/\La$} in the denominator. 

As the final step, we just need to redo Riemann's calculations in both 4.\ref{eq:e4} and 4.\ref{eq:e5} \textbf{---}  but now without forgetting any terms. In particular, when Riemann adds the series 4.\ref{eq:e2} to 4.\ref{eq:e1}, we should rather have \vspace{0.1mm}
\begin{align*}
\sum_{n=1}^{\La/2} \frac{1}{2n-1} - \sum^{\La/4}_{n=1} \frac{1}{4n-2} &= \sum_{n=1}^{\La/2} \frac{1}{2n-1} - \bigg( \sum^{\La/2}_{n=1} \frac{1}{4n-2} - \frac{\ln(2)}{4} \bigg) \\[5pt]
												      &= \sum^{\La/2}_{n=1} \frac{1}{4n-2} +  \frac{\ln(2)}{4}
\end{align*} 
And when Riemann adds 4.\ref{eq:e3} to this new series, we conclude instead with
\begin{align*}
\sum^{\La/2}_{n=1} \frac{1}{4n-2} +  \frac{\ln(2)}{4} -  \sum^{\La/4}_{n=1} \frac{1}{4n} &= \sum^{\La/2}_{n=1} \frac{1}{4n-2} +  \frac{\ln(2)}{4} - \bigg( \sum^{\La/2}_{n=1} \frac{1}{4n} -  \frac{\ln(2)}{4}  \bigg) \\[4.7pt]
														&= \frac{1}{2} \cdot \sum^{\La}_{n=1} \frac{(-1)^{n+1}}{n} + \frac{\ln(2)}{2} \\[5.5pt]
														&= \ln(2) 
\end{align*} 

There's no limit to the number of paradoxes that Riemann's theorem can generate. But, fortunately, we can always blame the same error (i.e. adding series of different length), and their resolutions are thus very similar. This is yet further evidence supporting the validity of our theory.   \\[12pt]

\setcounter{equation}{0}
\setcounter{section}{5}

\section*{5. Exponential }
\vspace{0.6mm}
Currently, in mathematics, there's a very strong consensus: the set {\small$\BB$}, which contains all infinite binary strings, is uncountable. According to Cantor's theory, the cardinality of {\small$\BB$} is so enormous that, paradoxically, it defies all counting methods \textbf{---} the very tools by which we can define a cardinality in the first place. Unsurprisingly perhaps, we'll use this last section to demonstrate the opposite. In particular, we'll set up a counting procedure that will reach all strings in {\small$\BB$}, and thus give a numerical value to its cardinality. Our conclusion will be simple: while {\small$\Na$} is too small to have a one-to-one correspondence with those strings, there's a larger set that can achieve this. 

To measure {\small$\BB$}'s cardinality, however, we must first specify the length of those strings \textbf{---} that is, how many digits they're made of. A natural choice is to take them with {\small$\La$} digits, one for each number in {\small$\Na$}. Admittedly, this is arbitrary, but our method is applicable to strings of any size. In fact, we just want to partition {\small$\BB$} into the smaller subsets {\small$\bb_n$}, which will contain, respectively, all strings with {\small$n$} ones and {\small$\La-n$} zeroes. It allows us to write 
\begin{equation} \label{eq:symmetry}
| \BB | = \sum_{i=0}^{\lambda}  |\bb_n| \\[6pt]
\end{equation}

The first term is simply {\small$|\bb_0|=1$}, since the string {\small{$000\ldots0_\sLa$} is the only one that contains no $1$. By just looking at the strings below, it's also easy to evaluate {\small$|\bb_1|$}.  \vspace{0.05mm}
\begin{align*}
10000000\ldots0_{\scriptscriptstyle\lambda} \\
01000000\ldots0_{\scriptscriptstyle\lambda} \\
00100000\ldots0_{\scriptscriptstyle\lambda} \\
\ldots \quad\quad\;\;\;\;\:  \\
00000000\ldots1_\lambda  
\end{align*} 
There's one different string for each digit where the $1$ can be, and there are {\small$\La$} digits. Consequently, {\small$| \bb_1| = \La$}. 

Now, for each string in {\small$ \bb_1$}, we'll count the digits where a second $1$ can be; this will give us {\small$| \bb_2|$}. To avoid counting twice the same string, however, this second $1$ can only access the digits to the right of the first $1$. Below, we're showing all the strings in {\small$ \bb_2$} that are made from the first three of {\small$ \bb_1$}. \vspace{0.48mm}
\begin{align*}
11000000\ldots0_{\scriptscriptstyle\lambda}   &&      01100000\ldots0_{\scriptscriptstyle\lambda}     &&  00110000\ldots0_{\scriptscriptstyle\lambda} \\
10100000\ldots0_{\scriptscriptstyle\lambda}   &&     01010000\ldots0_{\scriptscriptstyle\lambda}   &&  00101000\ldots0_{\scriptscriptstyle\lambda} \\
10010000\ldots0_{\scriptscriptstyle\lambda}    &&     01001000\ldots0_{\scriptscriptstyle\lambda}    &&  00100100\ldots0_{\scriptscriptstyle\lambda} \\
\ldots \quad\quad\;\;\;\;\: 				    &&			\ldots	\quad\quad\;\;\;\;\: 		&&  \ldots \quad\quad\;\;\;\;\:  \\
10000000\ldots1_{\scriptscriptstyle\lambda}    &&     01000000\ldots1_{\scriptscriptstyle\lambda}     && 00100000\ldots1_{\scriptscriptstyle\lambda} 
\end{align*}
The first block has its strings made from {\small{$100\ldots0_\sLa$}, and there are {\small$\La-1$} digits where a second $1$ can go; we therefore count that many strings. In the next block, the strings are made from {\small{$010\ldots0_\sLa$}, and we thus count {\small$\La-2$} more strings. Finally, in the last block, there are {\small$\La-3$} strings. Using the same counting method on all the other strings of {\small$\bb_1$}, we ultimately get 
 \begin{align}
 \begin{split}  \label{eq:f1}
 |\bb_2| &= \sum_{n=1}^{\lambda} (\lambda-n) \\[4.5pt]
 				      &= \lambda^2 - \bigg( \frac{ \lambda ^2 + \lambda }{2} \bigg) \\[7.5pt]
				      &= \frac{\lambda(\lambda-1)}{2} = \binom{\lambda}{2} \\[4.5pt]
 \end{split} 
 \end{align}
 Incidentally, since we're not using the {\small$+ \cdots$} notation, we can use the properties of {\small$\La$} to read the sum {\small$\sum_{n=1}^{\lambda} (\lambda-n)$} from left to right. This allows us to also write \vspace{0.5mm}
 \begin{equation}
 |\bb_2|  = \sum_{n=0}^{\La-1} n = 0 +1 + 2 + 3+ \cdots + (\La-1) = \binom{\lambda}{2} \\[2pt]
 \end{equation}
 
To continue further, we can now evaluate {\small$ |\bb_3| $} with the same approach: for each string in {\small$ \bb_1 $}, we simply count all the strings that can be made with two more $1$s. In particular, since there are {\small$\La-n$} available digits in the $n$th string of {\small$ \bb_1$}, we count {\small$\binom{\lambda-n}{2}$} strings from {\small$ \bb_3$}. Going through all the {\small$\La$} strings of {\small$ \bb_1$}, we obtain
\begin{align}
 \begin{split} 
 |\bb_3|  &= \sum_{n=1}^{\La} \binom{\La-n}{2}  \\[5pt]
		&=  \sum_{n=0}^{\La-1 }  \binom{n}{2}  = \binom{\La}{3} \\[4pt]
\end{split}
\end{align}
where the second line comes from reading the first sum from left to right; the terms are thus $0$, $0$, $1$, $3$, $6$, $10$, $15$, and so on.

Finally, by using an induction procedure with the same idea, we can generalize these results to {\small$|\bb_n| = \binom{\La}{n}$}. Nevertheless, for a full understanding of {\small$\BB$}'s structure, it's very important to notice that
\begin{equation} \label{general}
|\bb_n|  = \binom{\La}{n} =  \frac{\La!}{n!(\La-n)!}  = \binom{\La}{\La-n} = |\bb_{\La-n}| \\[2pt]
\end{equation}
This identity becomes evident after looking into the subset {\small$\bb_{\scriptscriptstyle\La-1}$}, where the strings have {\small$\La-1$} ones and only $1$ zero. These strings are \vspace{0.4mm}
\begin{align*}
\quad 011111111\ldots1_{\scriptscriptstyle\lambda} \\
\quad101111111\ldots1_{\scriptscriptstyle\lambda} \\
\quad110111111\ldots1_{\scriptscriptstyle\lambda} \\
\ldots \quad\quad\;\;\;\;\:  \\
\quad111111111\ldots0_\lambda 
\end{align*} 
The subsets {\small$\bb_1$ and {\small$\bb_{\La-1}$ are thus the same, except for the $0$s and the $1$s being interchanged \textbf{---} and that's also true for any pair {\small$\bb_n$ and {\small$\bb_{\La-n}$. The set {\small$\BB$} has therefore the same mirror-symmetry that's in the Pascal Triangle, where the left half is reflecting its right counterpart. In fact, the cardinality of {\small$\BB$} is equal to the sum of all its entries.
\begin{align} 
\begin{split} \label{eq:bsize}
\quad \quad \;\; |\BB|  &=  \sum_{n=0}^\La |\bb_n| \\
				  &= \: |\bb_0| + 
 		 \begin{pmatrix}
		\:\:\:\:\:1+1+1+1+1+\:1+\cdots+{\scriptscriptstyle\binom{\lambda}{0} } \: \\[2pt]
				  \:+\: 0+1+2+3+4+\:5+\cdots+{\scriptscriptstyle\binom{\lambda}{1} } \: \\[2pt]
				\:+\:0+0+1+3+6+10+\cdots+{\scriptscriptstyle\binom{\lambda}{2} }  \\[2pt]
				\:+\:0+0+0+1+4+10+\cdots+{\scriptscriptstyle\binom{\lambda}{3} }  \\[2pt]
				   \:\:  \cdots   \\[1.5pt]
				   +\: 0+0+0+0+0+0+\cdots+{\scriptscriptstyle\binom{\lambda}{\La} } 
  \end{pmatrix} \\[2pt]
  =& \sum_{n=0}^\La \binom{\La}{n}   = 2^\La \\[5pt]
  \end{split}
\end{align}

We must also notice the subset where each string has equally many $0$s and $1$s, namely {\small$\bb_{\La/2}$}. In the finite case, there are approximately {\small$\binom{n}{n/2} = \sqrt{2/n\pi} \cdot 2^n $} strings in this subset, and this becomes increasingly exact as $n$ grows bigger, due to the Stirling approximation. After taking the infinite limit, there's no reason for this numerical structure to suddenly disappear, and we therefore conclude with \vspace{0.8mm}
\begin{equation}
|\bb_{\La/2}| = \binom{\La}{\La/2} = \sqrt{\frac{2}{\La\pi}} \cdot 2^\La
\end{equation}
This is clearly consistent with the famous gaussian integral {\small$\int_{-\infty}^\infty e^{-x^2} dx = \sqrt{\pi}$}, which follows from the de Moivre–Laplace theorem. In particular, it explains that the normal distribution may be used as an approximation to the binomial distribution, which is essentially just counting binary strings. 

At this point, it's easy to expand our discussion to the set {\small$\BB_n$}, where the strings are made from $n$ different symbols. In fact, by using our previous results, we can quickly evaluate {\small$|\BB_3| = \sum_{m=0}^\La |\bb_m|$}, with the subsets {\small$ \bb_m$} now have strings with $m$ twos and {\small$\La-m$} of either zeros and ones. 

Since the twos occupy $m$ digits on a string, we use 5.\ref{general} to count {\small$\binom{\La}{m}$} combinations of them. For the {\small$\La-m$} remaining digits, however, we must refer to 5.\ref{eq:bsize}, and we thus count {\small$2^{\La-m}$} substrings that are made from zeroes and ones. The cardinality of {\small$\BB_3$} is therefore given as 
\begin{equation}
|\BB_3| = \sum_{m=0}^\La \binom{\La}{m} 2^{\La-m} = 3^\La \\[3.93pt]
\end{equation}
By using induction with the same counting method, it's straightforward to generalize our results to {\small$|\BB_n|=n^\La$}. The conclusion is thus that {\small$|\BB_i|<|\BB_j|$} for {\small$i<j$}. 

Cantor's theory would rather say that all those sets have the same cardinality, since they're supposedly all uncountable. But, by using {\small$\La$}, we've seen the complete opposite: {\small$\BB_2$} has a well-defined structure that can be used to correctly evaluate its cardinality \textbf{---} and, at this point, it's just intuitive that {\small$m^\La<n^\La$} for {\small$m<n$}. 

With that being said, before we're even allowed to conclude this work, we must discuss Cantor's famous diagonal argument. Although we disagree with most of his theory, there's at least one thing we do agree on: the set {\small$\BB$} is bigger than {\small$\Na$}, and the diagonal argument is a valid way to demonstrate this. The idea is to make a list of infinite binary strings, and first assume a one-to-one correspondence with the set {\small$\Na$}. From this assumption, it's possible to do the following enumeration.
\begin{align*}
1 \: & \leftrightarrow \: \underline{0}101011101001\ldots \\
2 \: & \leftrightarrow \: 1\underline{0}01001101001\ldots \\
3 \: & \leftrightarrow \: 01\underline{1}1010101101\ldots \\
4 \: & \leftrightarrow \: 101\underline{0}100100101\ldots \\
5 \: & \leftrightarrow \: 0010\underline{1}10001010\ldots \\
&\quad\quad\quad\quad  \ldots
\end{align*}
Now, if there's truly a one-to-one mapping between {\small$\Na$} and {\small$\BB$}, then the above list must contain all of {\small$\BB$}'s strings. However, by flipping every underlined digit, we make the string {\small$11010\ldots$}, and it differs by at least one digit from all the other ones. Cantor is thus right: this new string is not in the list, and {\small$\Na$} must be smaller than {\small$\BB$}. 

In fact, we can quickly achieve the same conclusion. Because of the one-to-one mapping with {\small$\Na$}, the above list contains {\small$\La$} strings \textbf{---} and since these are infinite, they each have a minimum of {\small$\La$} digits. But we can distinguish one string from another by only one different digit. Consequently, since the list contains {\small$\La^2$} digits to specify only {\small$\La$} strings, we must conclude that more strings exist outside of this list; it demonstrates that {\small$|\Na| < |\BB|$}. 

Cantor went further, however, and argued that only two types of infinity can exist. Either there's a one-to-one correspondence with the set {\small$\Na$}, and this is a \textit{countable} infinity; or there's no such mapping, and it's \textit{uncountable}. Cantor has therefore a binary view of infinity, and it simply fails to capture the great complexity that mathematics requires. Indeed, the complete opposite is true: there are infinitely many ways to use {\small$\La$} arithmetically, and they each are a different infinite quantity. 

There's even yet another problem in his theory. In particular, the ordinal number {\small$\omega^\omega$} \textbf{---} the equivalent of {\small$\La^\La$} in our work \textbf{---}  is said to have a one-to-one correspondence with the set {\small$\Na$}. But Cantor also understood intuitively that {\small$|\BB| = 2^\La$}, since he actually wrote {\small$|\BB| = 2^{\aleph_0}$} with {\small$|\Na| = \aleph_0$} in his theory. It's now clear that the inequality {\small$2^\La < \La^\La$} is true, however. Therefore, the set {\small$\BB$} would thus be uncountable, even if its cardinality is smaller than the countable {\small$\omega^\omega =\La^\La$}. 

Such major misconceptions always lead mathematicians to perceive infinity with a counterintuitive nature, and this surely explains why they came to accept so many paradoxical results. From 1873 to 1897, Cantor worked on developing  set theory and, thereafter, it rapidly became accepted as the foundation of mathematics \textbf{---} despite all its paradoxical results, and despite a strong resistance from many eminent mathematicians, such as Poincaré, Weyl, Brouwer and Kronecker \cite{problemshistory, Brouwer, weyl}. But, from our work, it's now evident that set theory is severely inconsistent with Calculus, and this cannot remain acceptable any longer. 

To fix the problems in Cantor's theory, we had to first recognize the main issue: the distinction made between cardinal and ordinal numbers. Nevertheless, by combining them into our own number {\small$\La$}, we could mostly recover the arithmetics of finite numbers. This allowed us to propose a theory of infinity that's compatible with Calculus, since it's no longer just countable or not. Incidentally, by working with {\small$\La$}, we were able to use infinity with high precision in our calculations, and this lead us to resolve several long-standing paradoxes in mathematics. Along the way, we therefore had to refute most results from set theory, which includes the alleged one-to-one mappings between the sets {\small$\Na$}, {\small$\ZZ$} and {\small$\QQ$}.

As our first next step, we'll shortly complete this work with a final section on {\small$\RR$} and the continuum, where {\small$|\BB_n|<|\RR|$} will be demonstrated. We'll thus give a definitive answer to the Continuum Hypothesis, since infinitely many quantities are in between the cardinalities of {\small$\Na$} and the real line \textbf{---}  a few examples are  {\small$\La+n$}, {\small$n\La$} and {\small$\La^n$}. Our goal with this theory, however, is to develop new mathematical tools, and ultimately make scientific discoveries. This will be the purpose of our upcoming papers. 

\medskip
\medskip
\medskip
\medskip
\medskip
\medskip
\medskip
\medskip

\bibliography{citations}{}
\bibliographystyle{abbrv}

\end{document}